\documentclass[12pt]{article}

\usepackage{amsfonts}

\usepackage{graphicx}

\newcommand{\ve}{\varepsilon}

\begin{document}

\begin{center}
 {\large \bf Spectral analysis for differential systems with a singularity.
 } \\[0.2cm]
 {\bf Mikhail Ignatyev} \\[0.2cm]
\end{center}

\thispagestyle{empty}

{\bf Abstract.}  We consider the differential system
$
y'-x^{-1}Ay-q(x)y=\rho By
$
with $n\times n$ matrices $A,B, q(x)$, where $A,B$ are constant, $B$ is diagonal, $A$ and $q(x)$ are off-diagonal, $q(\cdot)\in W^1_1[0,\infty)$. Some distinguished fundamental system of solutions is constructed.
Also, we discuss the inverse scattering problem and obtain the uniqueness result. \\
[0.1cm]



{\bf 1. Introduction.} We consider the differential system
$$
y'-x^{-1}Ay-q(x)y=\rho By, \ x>0, \eqno(1.1)
$$
with $n\times n$ matrices $A,B, q(x), x\in(0,\infty)$, where $A,B$ are constant.

Differential equations with coefficients having non-integrable singularities at the end or inside the interval often appear in various areas of natural sciences and engineering. For $n=2$, there exists extensive literature devoted to different aspects of spectral theory of the radial Dirac operators, see, for instance {\cite{TKBDir}}, {\cite{AHM}}, {\cite{AHM1}}, {\cite{Ser}}, {\cite{GoYu}}.

Systems of the form (1.1) with $n>2$ and arbitrary complex eigenvalues of the matrix $B$ appear to be considerably more difficult for investigation even in the "regular" case $A=0$ {\cite{BCsyst}}.  Some difficulties of principal matter also appear due to the presence of singularity at $x=0$. Whereas the "regular" case $A=0$ has been studied fairly completely to date {\cite{BCsyst}}, {\cite{Zh}}, {\cite{Ysyst}}, for system (1.1) with $A\neq 0$  there are no similar general results.

In this paper, we concentrate mostly on the construction of a distinguished basis of generalized eigenfunctions for (1.1). We call them  the {\it Weyl-type solutions}. The Weyl-type solutions play a central role in studying both direct and inverse spectral problems (see, for instance, {\cite{BeDT}},{\cite{Ybook}}). In presence of singularity at $x=0$ this step encounters some difficulties that do not appear in a "regular" case $A=0$. In particular, one can not use the auxiliary Cauchy problems with the initial conditions at $x=0$. The approach presented in {\cite{YIP93}} (see also {\cite{Fed}} and references therein) for
the scalar differential operators
$$
\ell y=y^{(n)}+\sum\limits_{j=0}^{n-2} \left(\frac{\nu_j}{x^{n-j}}+q_j(x)\right)y^{(j)}
$$
is based on some special solutions of the equation $\ell y=\lambda y$ that also satisfy certain Volterra integral equations. This approach assumes some additional decay condition for the coefficients $q_j(x)$ as $x\to 0$. In this paper, we do not impose any additional restrictions of such a type. Instead, we use the modification of the approach first presented in {\cite{BeDT}}
for the higher-order differential operators with {\it regular} coefficients on the whole line.

In brief outline our approach can be described as follows. We consider some auxiliary systems with respect to the functions with  values  in the exterior algebra $\wedge \mathbb{C}^n$. Our studying of these auxiliary systems centers on two families of  their solutions that also satisfy some asymptotical conditions as $x\to 0$ and $x\to\infty$ respectively, and can be constructed as solutions of certain {\it Volterra} integral equations. As in {\cite{BeDT}} we call these distinguished tensor solutions  the {\it fundamental tensors}.
The main difference from the above-mentioned method used in {\cite{YIP93}} is that we use the integral equations to construct the fundamental tensors rather than the solutions for the original system. Since each of the fundamental tensors has minimal growth (as $x\to 0$ or $x\to\infty$) among the solutions of the same auxiliary system,
  this step does not require any decay of $q(x)$ as $x\to 0$. As a next step, we show that the fundamental tensors are decomposable. Moreover, they can be represented as the wedge products of some solutions of original system (1.1) and these solutions can be shown to be the Weyl-type solutions of (1.1). This strategy allows us to construct the Weyl-type solutions via purely algebraic procedure and investigate their properties, in particular, their analytical properties with respect to the spectral parameter $\rho$ and asymptotical behavior as $\rho\to\infty$ and $\rho\to 0$.

The results obtained can be used in various areas of spectral theory. As an example, we consider the inverse scattering problem for (1.1). For the sake of brevity, we restrict our considerations here only with the uniqueness result, moreover, we assume that the discrete spectrum is empty. A general case and a constructive procedure for solving the inverse problem are planned to be considered in our further works.

Throughout the paper we assume the following

\medskip
{\bf Assumption 1.} $B=diag(b_1,\dots,b_n)$, $n>2$ , $A$ and $q(x)$ are off-diagonal, $q(\cdot)\in W^1_1[0,\infty)$ .
The eigenvalues $\{\mu_j\}_{j=1}^n$ of matrix $A$ are distinct and such that $\mu_j-\mu_k \notin \mathbb{Z}$ for $j\neq k$, moreover, $\mbox{Re}\mu_1<\mbox{Re}\mu_2<\dots<\mbox{Re}\mu_n$. The entries $b_1,\dots,b_n$ of matrix $B$ are nonzero distinct complex numbers such that  $\sum\limits_{j=1}^n b_j=0$.

\bigskip
{\bf 2. Solutions of the unperturbed system.}

Here we briefly discuss the unperturbed system:
$$
y'-x^{-1}Ay=\rho By \eqno(2.1)
$$
and introduce some fundamental systems of its solutions.

We start with considering (2.1) for $\rho=1$
$$
y'-x^{-1}Ay=By \eqno(2.2)
$$
but for {\it complex} values of $x$.

Let $\Sigma$ be the following union of lines through the origin in $\mathbb{C}$:
$$
\Sigma=\bigcup\limits_{(k,j): j\neq k}\left\{x:\mbox{Re}(xb_j)=\mbox{Re}(x b_k)\right\}.
$$
Then $\mathbb{C}\setminus\Sigma$ can be presented as a union of the sectors $S_\nu, \nu=\overline{1,N}$, where $N$ depends on the numbers $b_1,\dots,b_n$ and can be from 2 up to $n(n-1)$. For definiteness, we assume that $N>2$.

Consider some (arbitrary) sector $S_\nu$.It is well-known that there exists the ordering $R_1,\dots, R_n$ of the numbers $b_1,\dots,b_n$ such that $\mbox{Re}(R_1 x)<\mbox{Re}(R_2 x)\dots<\mbox{Re}(R_n x)$ for any $x\in S_\nu$. For $x\in\overline S_\nu\setminus\{0\}$ we consider the following fundamental matrices for system (2.2):
\begin{itemize}
\item $c(x)=(c_1(x),\dots,c_n(x))$, where
$$
c_k(x)=x^{\mu_k}\hat c_k(x),
$$
$\det c(x)\equiv 1$ and all $\hat c_k(\cdot)$ are entire functions, $\hat c_k(0)=\mathfrak{h}_k$, $\mathfrak{h}_k$ is an eigenvector of the matrix $A$ corresponding to the eigenvalue $\mu_k$;
\item $e(x)=(e_1(x),\dots,e_n(x))$, where
$$
e_k(x)=\mbox{e}^{\rho x R_k}(\mathfrak{f}_k+x^{-1}\eta_k(x)), \ \|\eta_k(x)\|\leq C_k, |x|\geq 1, x\in\overline S_\nu,
$$
$(\mathfrak{f}_1,\dots,\mathfrak{f}_n)=(\mathfrak{e}_1,\dots,\mathfrak{e}_n)\Pi_\nu$, $\{\mathfrak{e}_1,\dots,\mathfrak{e}_n\}$ is a standard coordinate column basis in $\mathbb{C}^n$ and $\Pi_\nu$ is a permutation matrix such that $(R_1,\dots,R_n)=(b_1,\dots,b_n)\Pi_\nu$.
\end{itemize}

\medskip
{\bf Condition $R_0$.} For all $\nu=\overline{1,N}$, $k=\overline{1,n}$ the numbers $$\Delta^0_k:=\det(e_1(x),\dots,e_{k-1}(x),c_k(x),\dots,c_n(x))$$ are not equal to 0.

\medskip
Under Condition $R_0$ for $x\in\overline S_\nu\setminus\{0\}$ there exists (and is unique) the fundamental matrix $\psi^0(x)$ such that
$$
\psi^0_k(x)=\mbox{e}^{\rho x R_k}(\mathfrak{f}_k+o(1)), x\to\infty, x\in S_\nu, \ \psi^0_k(x)=O(x^{\mu_k}), x\to 0.
$$

\medskip
{\bf Remark 2.1.} The fundamental matrix $c(x)$ is uniquely determined by conditions specified above provided that we fix some set of eigenvectors $\{\mathfrak{h}_k\}_{k=\overline{1,n}}$ of matrix $A$ and some branch for $\arg x$ in $S_\nu$. The fundamental matrix $e(\cdot)$ is uniquely determined in any {\it Stokes sector} containing $\overline S_\nu\setminus \{0\}$ but not in $S_\nu$ itself. But this does not affect neither Condition $R_0$ nor our definition of $\psi^0(\cdot)$. Indeed, $e(x)$ for $x\in S_\nu$ is determined uniquely up to multiplication by upper-triangular constant matrix with diagonal entries equal to 1.It is clear that such a multiplication does not change the tensors $e_1(x)\wedge\dots\wedge e_{k-1}(x)$ and therefore the numbers $\Delta^0_k$. Furthermore, for $\psi^0(x)$ one has the representation $\psi^0(x)=c(x)l$, where constant lower-triangular matrix $l$ arises from the $LU$-factorization of the matrix $c^{-1}(x)e(x)$. It is clear that multiplication of $e(x)$ by upper-triangular constant matrix with diagonal entries equal to 1 does not change the matrix $l$.

\medskip
Now we return to (2.1) with arbitrary $\rho\in \overline S_\nu\setminus\{0\}$ and real positive $x$. Notice that if $y(x)$ satisfies (2.2) then $Y(x,\rho):=y(\rho x)$ satisfies (2.1). Taking this into account, we define the matrix solutions $C(x,\rho)$, $E(x,\rho)$, $\Psi^0(x,\rho)$ of (2.1) as follows: $C(x,\rho):=c(\rho x)$, $\Psi^0(x,\rho):=\psi^0(\rho x)$, $E(x,\rho):=e(\rho x)$ .

\bigskip
{\bf 3. Fundamental tensors.}

In this section we consider the following auxiliary equations:
$$
Y'=Q^{(m)}(x,\rho)Y, \eqno(3.1)
$$
where $Y$ is a function with values in the exterior product $\wedge^m \mathbb{C}^n$. Here and below
$$
Q(x,\rho):=x^{-1}A+\rho B+q(x)
$$
and for given $n\times n$ matrix $M$  $M^{(m)}$ denotes the operator acting in $\wedge^m \mathbb{C}^n$ so that for any vectors $u_1, \dots, u_m$ the following identity holds {\cite{BeDT}}:
$$
M^{(m)}(u_1\wedge u_2 \wedge\dots\wedge u_m)=\sum\limits_{j=1}^m u_1\wedge u_2 \wedge\dots\wedge u_{j-1}\wedge M u_j\wedge u_{j+1}\wedge\dots\wedge u_m.
$$
In what follows, we shall also use the following notations.

We denote by $\mathcal{A}_{m}$ the set of all ordered multi-indices $\alpha=(\alpha_1, \dots, \alpha_m)$, $\alpha_1<\alpha_2<\dots<\alpha_m$, $\alpha_j\in\{1,2,\dots,n\}$. For a set of vectors $u_1, \dots, u_n$
 from $\mathbb{C}^n$ and a multi-index $\alpha\in \mathcal{A}_m$ we define
 $$u_\alpha:=u_{\alpha_1}\wedge\dots\wedge u_{\alpha_m}.$$
Let $a_1,\dots,a_n$ be a numerical sequence. For $\alpha\in\mathcal{A}_m$ we define $$a_\alpha:=\sum\limits_{j\in\alpha} a_j,$$ for $k\in\overline{1,n}$ we denote $$\overrightarrow{a}_k:=\sum\limits_{j=1}^k a_j,\ \overleftarrow{a}_k:=\sum\limits_{j=k}^n a_j.$$
For a multi-index $\alpha$ the symbol $\alpha'$ denotes the ordered multi-index that complements $\alpha$ to $(1,2,\dots,n)$.
We note that Assumption 1 implies, in particular, that $\sum\limits_{k=1}^n \mu_k=\sum\limits_{k=1}^n R_k=0$ and therefore for any  multi-index $\alpha$ one has $R_{\alpha'}=-R_\alpha$ and $\mu_{\alpha'}=-\mu_\alpha$.

It is clear that if a set $\{u_1, \dots, u_n\}$ is a basis
 in $\mathbb{C}^n$ then the set $\{u_\alpha\}_{\alpha\in\mathcal{A}_m}$ is a basis in $\wedge^m \mathbb{C}^n$. In particular the set $\{\mathfrak{e}_\alpha\}_{\alpha\in\mathcal{A}_m}$, where $\{\mathfrak{e}_k\}_{k=1}^n$ is a standard basis in $\mathbb{C}^n$ is a basis in $\wedge^m \mathbb{C}^n$. In $\wedge^m \mathbb{C}^n$ we define a norm as follows:
$$
\left\|\sum\limits_{\alpha\in\mathcal{A}_m}h_\alpha \mathfrak{e}_\alpha\right\|:=\sum\limits_{\alpha\in\mathcal{A}_m}|h_\alpha| .
$$
Then for $n\times n$ matrix $M$ one has:
$$
\|M^{(m)}\|\leq m\|M\|. \eqno(3.2)
$$
For $h\in \wedge^n \mathbb{C}^n$ we define $|h|$ as a constant in the following representation:
$$
h=|h|\mathfrak{e}_1\wedge\mathfrak{e}_2\wedge\dots\wedge\mathfrak{e}_n.
$$
In the sequel we assume that $\rho\in\overline S_\nu$ for some arbitrary fixed $\nu$. We consider the following Volterra equations:
$$
Y(x)=T_k^0(x,\rho)+\int\limits_0^x G_{n-k+1}(x,t,\rho)\left(q^{(n-k+1)}(t)Y(t)\right)dt, \eqno(3.3)
$$
$$
Y(x)=F_k^0(x,\rho)-\int\limits_x^\infty G_{k}(x,t,\rho)\left(q^{(k)}(t)Y(t)\right)dt, \eqno(3.4)
$$
where
$$T_k^0(x,\rho):=C_k(x,\rho)\wedge\dots\wedge C_n(x,\rho),$$
$$F_k^0(x,\rho):=E_1(x,\rho)\wedge\dots\wedge E_k(x,\rho)=\Psi^0_1(x,\rho)\wedge\dots\wedge \Psi^0_k(x,\rho)$$
and $G_m(x,t,\rho)$ is operator acting in $\wedge^m \mathbb{C}^n$ as follows:
$$
G_m(x,t,\rho)f=
\sum\limits_{\alpha\in\mathcal{A}_m}(-1)^{\sigma_\alpha}\left|f\wedge\Psi^0_{\alpha'}(t,\rho)\right|\Psi^0_{\alpha}(x,\rho)=
\sum\limits_{\alpha\in\mathcal{A}_m}(-1)^{\sigma_\alpha}\left|f\wedge E_{\alpha'}(t,\rho)\right|E_{\alpha}(x,\rho).
$$
Here $\sigma_\alpha\in\{0,1\}$ is such that $(-1)^{\sigma_\alpha}=\left|\mathfrak{f}_\alpha\wedge\mathfrak{f}_{\alpha'}\right|$.

It follows directly from the definition of the fundamental matrix $\Psi^0(x,\rho)$ (see the previous section) that for any $\alpha\in\mathcal{A}_m$ the following estimates hold:
$$
\|\Psi^0_\alpha(x,\rho)\|\leq C \left\{ \begin{array}{l}
\left|(\rho x)^{\mu_\alpha}\right|, \ |\rho x|\leq 1 \\
\left|\exp(\rho x R_\alpha)\right|, \ |\rho x|>1.
                                                 \end{array}\right. \eqno(3.5)
                                                 $$
Here and below we use the same symbol $C$ to denote different constants that does not depend on $x,\rho$.

\medskip
{\bf Theorem 3.1.} For any $\rho\in\overline S_\nu\setminus\{0\}$ equations (3.3) and (3.4) have unique solutions $T_k(x,\rho)$ and $F_k(x,\rho)$ respectively such that
$$
\|T_k(x,\rho)\|\leq C \left\{ \begin{array}{l}
\left|(\rho x)^{\overleftarrow{\mu}_k}\right|, \ |\rho x|\leq 1 \\
\left|\exp(\rho x \overleftarrow{R}_k)\right|, \ |\rho x|>1.
\end{array}\right.
$$
$$
\|F_k(x,\rho)\|\leq C \left\{ \begin{array}{l}
\left|(\rho x)^{\overrightarrow{\mu}_k}\right|, \ |\rho x|\leq 1 \\
\left|\exp(\rho x \overrightarrow{R}_k)\right|. \ |\rho x|>1.
\end{array}\right.
$$
The following asymptotics hold:
$$
F_k(x,\rho)=\exp(\rho x \overrightarrow{R}_k)\left(\mathfrak{f}_1\wedge\dots\wedge\mathfrak{f}_k+o(1)\right), \ x\to\infty,
$$
$$
T_k(x,\rho)=(\rho x)^{\overleftarrow{\mu}_k}\left(\mathfrak{h}_k\wedge\dots\wedge\mathfrak{h}_n+o(1)\right), x\to 0.
$$

\medskip
{\bf Remark 3.1.} Notice that $G_m(x,t,\rho)$ is the Green's function for the nonhomogeneous  equation
$$
Y'=Q_0^{(m)}(x,\rho)Y+F, \ Q_0(x,\rho)=x^{-1}A+\rho B, \eqno(3.6)
$$
i.e., for any $x_0\in[0,\infty]$ the integral $Y(x)=\int\limits_{x_0}^x G_m(x,t,\rho)F(t)dt$ (if exists) solves (3.6). Thus the tensors $T_k(x,\rho)$, $F_k(x,\rho)$ are solutions for (3.1) with $m=n-k+1$ and $m=k$ respectively.

\medskip
{\bf Proof of Theorem 3.1.} Let us notice first that $T^0_k(x,\rho)$ admits the estimates:
$$
\|T^0_k(x,\rho)\|\leq C \left\{ \begin{array}{l}
\left|(\rho x)^{\overleftarrow{\mu}_k}\right|, \ |\rho x|\leq 1 \\
\left|\exp(\rho x \overleftarrow{R}_k)\right|, \ |\rho x|>1.
\end{array}\right. \eqno(3.7)
$$
Indeed, the first estimate follows directly from the definition of $T^0_k$, the second one can be obtained from the representation:
$$
T^0_k(x,\rho)=\sum\limits_{\alpha\in\mathcal{A}_{n-k+1}}T^0_{k\alpha} E_\alpha(x,\rho),
$$
where the numbers $T^0_{k\alpha}$ do not depend on $x,\rho$ and the observation that $\|E_\alpha(x,\rho)\|<C \left|\exp(\rho x \overleftarrow{R}_k)\right|$ for any $\alpha\in\mathcal{A}_{n-k+1}$, $|\rho x|>1$.

Now we apply the successive approximation method to (3.3) and consider the sequence $T^r_k(x,\rho)$ defined in recurrent way as follows:
$$
T^{r+1}_k(x,\rho)=\int\limits_0^x G_{n-k+1}(x,t,\rho)\left(q^{(n-k+1)}(t)T^r_k(t,\rho)\right)dt. \eqno(3.8)
$$
We use induction in $r$ to obtain the estimate:
$$
\|T^r_k(x,\rho)\|\leq C_0\cdot\frac{C^r}{r!}\cdot\left(\int\limits_0^x \|q(t)\|dt \right)^r \cdot \left\{ \begin{array}{l}
\left|(\rho x)^{\overleftarrow{\mu}_k}\right|, \ |\rho x|\leq 1 \\
\left|\exp(\rho x \overleftarrow{R}_k)\right|, \ |\rho x|>1.
\end{array}\right. \eqno(3.9)
$$
The base is obvious since (3.9) for $r=0$ is exactly estimate (3.7). Now let (3.9) be true for some $r$. Then for $|\rho x|\leq 1$ we have:
$$
\|T^{r+1}_k(x,\rho)\|\leq C_0\cdot\frac{C^r}{r!} \sum\limits_{\alpha\in\mathcal{A}_{n-k+1}} C_\alpha \int\limits_{0}^{x}
\left|\left(\frac{x}{t}\right)^{\mu_\alpha}\right|\cdot \left|(\rho t)^{\overleftarrow{\mu}_k}\right|\cdot\|q(t)\|\cdot
\left(\int_0^t \|q(\tau)\|d\tau\right)^r dt.
$$
Taking into account that for any $\alpha\in\mathcal{A}_{n-k+1}$ one has $\mbox{Re}(\overleftarrow{\mu}_k-\mu_\alpha)\geq 0$ and therefore $\left|(t/x)^{\overleftarrow{\mu}_k-\mu_\alpha}\right|\leq 1$ for $t\leq x$ we arrive at the first of desired estimates. Further, for $|\rho x|>1$ we write:
$$
\|T^{r+1}_k(x,\rho)\|\leq C_0\cdot\frac{C^r}{r!} \sum\limits_{\alpha\in\mathcal{A}_{n-k+1}} C_\alpha \int\limits_{0}^{|\rho|^{-1}}
|\exp(\rho x R_\alpha)|\cdot \left|(\rho t)^{\mu_{\alpha'}}\right|\cdot \left|(\rho t)^{\overleftarrow{\mu}_k}\right|\cdot\|q(t)\|\cdot
\left(\int_0^t \|q(\tau)\|d\tau\right)^r dt
$$
$$
+C_0\cdot\frac{C^r}{r!} \sum\limits_{\alpha\in\mathcal{A}_{n-k+1}} C_\alpha \int\limits_{|\rho|^{-1}}^x
|\exp(\rho x R_\alpha)|\cdot |\exp(\rho t R_{\alpha'})| \cdot |\exp(\rho t \overleftarrow{R}_k)| \cdot\|q(t)\|\cdot
\left(\int_0^t \|q(\tau)\|d\tau\right)^r dt \leq
$$
$$
\leq C_0\cdot\frac{C^r}{r!}\cdot|\exp(\rho x \overleftarrow{R}_k)| \sum\limits_{\alpha\in\mathcal{A}_{n-k+1}} C_\alpha \int\limits_{0}^{|\rho|^{-1}}
|\exp(\rho x (R_\alpha-\overleftarrow{R}_k))|\cdot \|q(t)\|\cdot
\left(\int_0^t \|q(\tau)\|d\tau\right)^r dt
$$
$$
+C_0\cdot\frac{C^r}{r!}\cdot|\exp(\rho x \overleftarrow{R}_k)| \sum\limits_{\alpha\in\mathcal{A}_{n-k+1}} C_\alpha \int\limits_{|\rho|^{-1}}^x
|\exp(\rho (x-t) (R_\alpha-\overleftarrow{R}_k))| \cdot\|q(t)\|\cdot
\left(\int_0^t \|q(\tau)\|d\tau\right)^r dt,
$$
where we have taken into account that as above $\left|(\rho t)^{\overleftarrow{\mu}_k-\mu_\alpha}\right|\leq 1$ if $|\rho t|\leq 1$. Observe that $\mbox{Re}(\rho(R_\alpha-\overleftarrow{R}_k))\leq 0$ and therefore in the estimate above one has $|\exp(\rho x (R_\alpha-\overleftarrow{R}_k))|\leq 1$, $|\exp(\rho (x-t) (R_\alpha-\overleftarrow{R}_k))|\leq 1$. Taking this into account, we obtain the desired estimate for $T^{r+1}_k(x,\rho)$, $|\rho x|>1$.

Estimate (3.9) guarantees a convergency of the successive approximation method and yields the desired estimate for $T_k(x,\rho)$. Using this estimate and equation (3.3) we obtain
 specified in theorem
 asymptotics for $T_k(x,\rho)$ as $x\to 0$.

This completes the proof in what concerns equation (3.3). Equation (3.4) can be considered in a similar way.
$\hfil\Box$

\medskip
In the rest of the section we consider the asymptotical behavior of fundamental tensors $T_k(x,\rho)$, $F_k(x,\rho)$ for $\rho\to\infty$ and $\rho\to 0$.

\medskip
{\bf Lemma 3.1.} Let  $T\in(0,\infty)$ be arbitrary fixed. Then for the function:
$$
H_0(x,\rho):=\int\limits_0^{|\rho|^{-1}}G_{n-k+1}(x,t,\rho)\left(q^{(n-k+1)}(t)T_k(t,\rho)\right)dt
$$
the following estimate holds:
$$
\|H_0(x,\rho)\|\leq C|\rho|^{-1}|\exp(\rho x\overleftarrow{R}_k)|
$$
for $|\rho|>T^{-1}$ uniformly in $x\in[|\rho|^{-1}, T]$, where the constant $C$ depends only on $T$.

\medskip
{\bf Proof.} Using the estimates for $T_k(x,\rho)$ from Theorem 3.1  we obtain for $|\rho x|\geq 1$:
$$
\|H_0(x,\rho)\|\leq C \sum\limits_{\alpha\in\mathcal{A}_{n-k+1}} \int\limits_0^{|\rho|^{-1}}
|\exp(\rho x R_\alpha)|\cdot \left|(\rho t)^{\mu_{\alpha'}}\right|\cdot\left|(\rho t)^{\overleftarrow{\mu}_k}\right|\cdot\|q(t)\|dt \leq
$$
$$
\leq C |\exp(\rho x \overleftarrow{R}_k)| \sum\limits_{\alpha\in\mathcal{A}_{n-k+1}} \int\limits_0^{|\rho|^{-1}}
|\exp(\rho x (R_\alpha-\overleftarrow{R}_k))|\cdot \left|(\rho t)^{\overleftarrow{\mu}_k-\mu_\alpha}\right|\cdot\|q(t)\|dt.
$$
As in proof of Theorem 3.1 we notice that $|\exp(\rho x (R_\alpha-\overleftarrow{R}_k))|\cdot \left|(\rho t)^{\overleftarrow{\mu}_k-\mu_\alpha}\right|\leq 1$ and therefore:
$$
\|H_0(x,\rho)\|\leq C \cdot |\exp(\rho x \overleftarrow{R}_k)| \cdot \max\limits_{x\in[0,|\rho|^{-1}]}  \|q(x)\| \cdot \sum\limits_{\alpha\in\mathcal{A}_{n-k+1}} \int\limits_0^{|\rho|^{-1}}
dt
$$
for $x\in[|\rho|^{-1},T]$.
$\hfil\Box$

\medskip
{\bf Lemma 3.2.} Let $T\in(0,\infty)$ be arbitrary fixed. Then:
\begin{enumerate}
  \item for any two multi-indices $\alpha,\beta\in\mathcal{A}_{n-k+1}$ the function:
$$
H^0_{\alpha\beta}(x,\rho)=\int\limits_{|\rho|^{-1}}^x
\left|\left(q^{(n-k+1)}(t)E_\alpha(t,\rho)\right)\wedge E_{\beta'}(t,\rho)\right| E_\beta(x,\rho) dt
$$
admits the estimate:
$$
\|H^0_{\alpha\beta}(x,\rho)\|\leq C|\rho^{-\varepsilon}\exp(\rho x\overleftarrow{R}_k)|
$$
for $|\rho|>T^{-1}$ uniformly in $x\in[|\rho|^{-1}, T]$, where $\varepsilon\in(0,1)$ is arbitrary and constant $C$ depends only on $\varepsilon$ and $T$;\\
  \item for any multi-index $\beta\in\mathcal{A}_{k}$ the function:
$$
H^\infty_{\beta}(x,\rho)=\int\limits_{|\rho|^{-1}}^x
\left|\left(q^{(k)}(t)F^0_k(t,\rho)\right)\wedge E_{\beta'}(t,\rho)\right| E_\beta(x,\rho) dt
$$
admits the estimate:
$$
\|H^\infty_{\beta}(x,\rho)\|\leq C|\rho^{-1}\exp(\rho x\overrightarrow{R}_k)|
$$
for $|\rho|>T^{-1}$ uniformly in $x\in[T,\infty)$, constant $C$ depends only on $T$.
\end{enumerate}
.

\medskip
{\bf Proof.} For any multi-index $\alpha\in\mathcal{A}_{m}$ $E_\alpha(x,\rho)$ admits the representation:
$$
E_\alpha(x,\rho)=\exp(\rho x R_\alpha)\left(\mathfrak{f}_\alpha+(\rho x)^{-1}\eta_\alpha(\rho x)\right),
$$
where $\|\eta_\alpha(z)\|\leq C_\alpha$ for any $z\in\overline S_\nu$, $|z|\geq 1$. Thus we can split $H^0_{\alpha\beta}(x,\rho)$ as follows:
$$
H^0_{\alpha\beta}(x,\rho)=h_0(x,\rho)+h_1(x,\rho),
$$
$$
h_1(x,\rho)=\int\limits_{|\rho|^{-1}}^x (\rho t)^{-1}\exp(\rho t R_\alpha)\cdot
\left|\left(q^{(n-k+1)}(t)\eta_\alpha(t,\rho)\right)\wedge E_{\beta'}(t,\rho)\right| E_\beta(x,\rho) dt,
$$
$$
h_0(x,\rho)=h_{00}(x,\rho)+h_{01}(x,\rho),
$$
$$
h_{01}(x,\rho)=\int\limits_{|\rho|^{-1}}^x (\rho t)^{-1}\exp(\rho t R_\alpha- \rho t R_\beta)\cdot
\left|\left(q^{(n-k+1)}(t)\mathfrak{f}_\alpha\right)\wedge \eta_{\beta'}(t,\rho)\right| E_\beta(x,\rho) dt,
$$
$$
h_{00}(x,\rho)=h_{000}(x,\rho)+h_{001}(x,\rho),
$$
$$
h_{001}(x,\rho)=\int\limits_{|\rho|^{-1}}^x (\rho x)^{-1}\exp(\rho t R_\alpha- \rho t R_\beta+\rho x R_\beta)\cdot
\left|\left(q^{(n-k+1)}(t)\mathfrak{f}_\alpha\right)\wedge \mathfrak{f}_{\beta'}\right| \eta_\beta(x,\rho) dt,
$$
$$
h_{000}(x,\rho)=\int\limits_{|\rho|^{-1}}^x \exp(\rho t R_\alpha- \rho t R_\beta+\rho x R_\beta)\cdot
\left|\left(q^{(n-k+1)}(t)\mathfrak{f}_\alpha\right)\wedge \mathfrak{f}_{\beta'}\right| \mathfrak{f}_\beta dt.
$$
We notice that for $t>|\rho|^{-1}$ and $\ve\in(0,1)$ one has $|\rho t|^{\ve-1}< 1$. Thus one can estimate:
$$
\|\rho^\ve\exp(-\rho x \overleftarrow{R}_k)h_1(x,\rho)\|\leq C
\int\limits_{|\rho|^{-1}}^x t^{-\ve} \|q(t)\|\left|\exp(\rho t R_\alpha-\rho x \overleftarrow{R}_k+\rho(x-t)R_\beta)\right|dt.
$$
 Let us transform: $t R_\alpha-x\overleftarrow{R}_k+(x-t)R_\beta=(x-t)(R_\beta-\overleftarrow{R}_k)+t(R_\alpha-\overleftarrow{R}_k)$ .Then taking into account that for $\rho\in\overline S_\nu$ and any $\alpha,\beta\in\mathcal{A}_{n-k+1}$ one has $\mbox{Re}(\rho(R_\alpha-\overleftarrow{R}_k))\leq 0$, $\mbox{Re}(\rho(R_\beta-\overleftarrow{R}_k))\leq 0$ we conclude that for $0<t<x$ we have $\left|\exp(\rho t R_\alpha-\rho x \overleftarrow{R}_k+\rho(x-t)R_\beta)\right|\leq 1$. Thus one gets:
$$
\|\rho^\ve\exp(-\rho x \overleftarrow{R}_k)h_1(x,\rho)\|\leq C \max\limits_{t\in[0,T]}\|q(t)\|\cdot
\int\limits_{0}^T t^{-\ve}dt.
$$
The functions $h_{01}(x,\rho)$ and $h_{001}(x,\rho)$ can be estimated in a similar way (in the last case one should notice first that $x^{-1}\leq t^{-1}$ for $t\in[|\rho|^{-1},x]$) and we obtain:
$$
\|H^0_{\alpha\beta}(x,\rho)-h_{000}(x,\rho)\|\leq C |\rho|^{-\varepsilon}|\exp(\rho x\overleftarrow{R}_k)|. \eqno(3.10)
$$
Now we consider $h_{000}(x,\rho)$ and we notice first that $h_{000}(x,\rho)\equiv 0$ if $\alpha=\beta$. Indeed, we have:
$$
q^{(n-k+1)}(t)\mathfrak{f}_\alpha=\sum\limits_{\gamma\in\mathcal{A}_{n-k+1}} q^{(n-k+1)}_{\gamma\alpha}(t)\mathfrak{f}_\gamma
$$
with complex (scalar) coefficients $q^{(n-k+1)}_{\gamma\alpha}(t)$. Therefore we can write:
$$
\left|\left(q^{(n-k+1)}(t)\mathfrak{f}_\alpha\right)\wedge\mathfrak{f}_{\beta'}\right|=\sum\limits_{\gamma\in\mathcal{A}_{n-k+1}} q^{(n-k+1)}_{\gamma\alpha}(t)\left|\mathfrak{f}_\gamma\wedge\mathfrak{f}_{\beta'}\right|. \eqno(3.11)
$$
If $\alpha=\beta$ (3.11) becomes:
$$
\left|\left(q^{(n-k+1)}(t)\mathfrak{f}_\alpha\right)\wedge\mathfrak{f}_{\alpha'}\right|=\sum\limits_{\gamma\in\mathcal{A}_{n-k+1}} q^{(n-k+1)}_{\gamma\alpha}(t)\left|\mathfrak{f}_\gamma\wedge\mathfrak{f}_{\alpha'}\right|=
q^{(n-k+1)}_{\alpha\alpha}(t)\left|\mathfrak{f}_\alpha\wedge\mathfrak{f}_{\alpha'}\right|=0,
$$
since $q_{jj}(t)\equiv 0, j=\overline{1,n}$ implies $q^{(m)}_{\alpha\alpha}(t)\equiv 0$ for all $m=\overline{1,n}$, $\alpha\in\mathcal{A}_m$.

Now, suppose $\alpha\neq \beta$. Then (3.11) yields:
$$
\|h_{000}(x,\rho)\|=\left\|\int\limits_{|\rho|^{-1}}^x \exp(\rho t R_\alpha- \rho t R_\beta+\rho x R_\beta)\cdot
q^{(n-k+1)}_{\beta\alpha}(t)\mathfrak{f}_\beta dt\right\|.
$$
Using integration by parts we obtain:
$$
\|h_{000}(x,\rho)\|\leq C |\rho|^{-1}\|q(x)\|\left|\exp (\rho x R_\alpha)\right|+C_1 |\rho|^{-1}\left|\exp (\rho x R_\beta)\right|+$$$$
C_2 |\rho|^{-1}\int\limits_{|\rho|^{-1}}^x \left|\exp(\rho t R_\alpha- \rho t R_\beta+\rho x R_\beta)\right|\cdot
|p^{(n-k+1)}_{\beta\alpha}(t)|dt,
$$
where $p^{(n-k+1)}_{\beta\alpha}(t)\in L(0,\infty)$. Therefore we have:
$$
\|\rho \exp(-\rho x \overleftarrow{R}_k) h_{000}(x,\rho)\|\leq
C_1 \left|\exp(\rho x (R_\alpha-\overleftarrow{R}_k))\right|+C_2 \left|\exp(\rho x (R_\beta-\overleftarrow{R}_k))\right|+$$$$
C_3 \int\limits_{|\rho|^{-1}}^x \left|\exp(\rho t R_\alpha- \rho t R_\beta+\rho x R_\beta-\rho x \overleftarrow{R_k})\right|\cdot
|p^{(n-k+1)}_{\beta\alpha}(t)|dt \leq C
$$
for $|\rho|>T^{-1}$, $x\in[|\rho|^{-1}, T]$ and $C$ depend only on $T$. Together with (3.10) this completes the proof of  the first part of the Lemma.

For the function $H^\infty_{\beta}(x,\rho)$ one can use the similar arguments. The only difference is that now we don't need to split $(\rho t)^{-1}$ as $(\rho t)^{-\varepsilon}\cdot(\rho t)^{\varepsilon-1}$. For instance, we can estimate (where $\alpha=(1,\dots,k)$):
$$
\|h^\infty_1(x,\rho)\|:=\left\|\int\limits_{x}^\infty (\rho t)^{-1}\exp(\rho t R_\alpha)\cdot
\left|\left(q^{(k)}(t)\eta_\alpha(t,\rho)\right)\wedge E_{\beta'}(t,\rho)\right| E_\beta(x,\rho) dt\right\| \leq$$$$
\leq C \left|\exp(\rho x \overrightarrow{R}_k)\right|
\int\limits_{x}^\infty |\rho t|^{-1} \|q(t)\|\left|\exp(\rho (t-x) \overrightarrow{R}_k+\rho(x-t)R_\beta)\right|dt \leq
$$$$
\leq C T^{-1}\left|\rho^{-1}\exp(\rho x \overrightarrow{R}_k)\right|
\int\limits_{T}^\infty  \|q(t)\|dt .
$$
$\hfil\Box$

\medskip
{\bf Lemma 3.3.} Let $T\in(0,\infty)$ be arbitrary fixed. Then:
\begin{enumerate}
  \item the following estimate holds:
  $$
\|T_k(x,\rho)-T_k^0(x,\rho)\|<C\left|\rho^{-\varepsilon}\exp(\rho x \overleftarrow{R}_k)\right| $$
  for $|\rho|>T^{-1}$ uniformly in $x\in[|\rho|^{-1}, T]$, where $\varepsilon\in(0,1)$ is arbitrary and the constant $C$ depend only on $\varepsilon$ and $T$;\\
  \item the following estimate holds:
$$
\|F_k(x,\rho)-F_k^0(x,\rho)\|<C\left|\rho^{-1}\exp(\rho x \overrightarrow{R}_k)\right|
$$
 for $|\rho|>T^{-1}$ uniformly in $x\in[T,\infty)$; the constant $C$ depend only on $T$.
\end{enumerate}

\medskip
{\bf Proof.} Let us fix an arbitrary $T>0$ and rewrite equation (3.3) as follows:
$$
\hat T_k(x,\rho)=\int\limits_{|\rho|^{-1}}^x G_{n-k+1}(x,t,\rho)\left(q^{(n-k+1)}(t)\hat T_k(t,\rho)\right)dt+\hat T^0_k(x,\rho), \eqno(3.12)
$$
where $\hat T_k(x,\rho):= T_k(x,\rho)- T^0_k(x,\rho)$, $\rho\in \overline S_\nu$ is arbitrary such that $|\rho|>T^{-1}$, $x\in[|\rho|^{-1}, T]$ and
$$
\hat T^0_k(x,\rho)=H_0(x,\rho)+H_1(x,\rho),
$$
$$
H_0(x,\rho)=\int\limits_0^{|\rho|^{-1}}G_{n-k+1}(x,t,\rho)\left(q^{(n-k+1)}(t)T_k(t,\rho)\right)dt,
$$
$$
H_1(x,\rho)=\int\limits_{|\rho|^{-1}}^x G_{n-k+1}(x,t,\rho)\left(q^{(n-k+1)}(t)T^0_k(t,\rho)\right)dt.
$$
We use the representations:
$$
G_{n-k+1}(x,t,\rho)f=\sum\limits_{\beta\in\mathcal{A}_{n-k+1}}(-1)^{\sigma_\beta}\left|f\wedge E_{\beta'}(t,\rho)\right|\cdot E_\beta(x,\rho),
$$
$$
T^0_k(t,\rho)=\sum\limits_{\alpha\in\mathcal{A}_{n-k+1}}T^0_{k\alpha} E_\alpha(t,\rho),
$$
where the numbers $T^0_{k\alpha}$ do not depend on $t,\rho$ and obtain:
$$
H_1(x,\rho)=\sum\limits_{\alpha,\beta\in\mathcal{A}_{n-k+1}} (-1)^{\chi_{k\alpha\beta}} T^0_{k\alpha} H^0_{\alpha\beta}(x,\rho), \ \chi_{k\alpha\beta}\in\{0,1\}.
$$
Lemma 3.2 yields the estimate:
$$
\|H_1(x,\rho)\|\leq C|\rho|^{-\ve}\left|\exp(\rho x \overleftarrow{R}_k)\right|.
$$
Together with the estimate for $H_0(x,\rho)$ from Lemma 3.1 this yields the estimate:
$$
\|\hat T^0_k(x,\rho)\|\leq C_0|\rho|^{-\ve}\left|\exp(\rho x \overleftarrow{R}_k)\right|,
$$
where $C_0$ depend only on $T$ and $\ve$ and do not depend on $x,\rho$. Then using the successive approximation method to solve (3.12) and proceeding as in the proof of Theorem 3.1 we obtain the estimate:
$$
\|\hat T_k(x,\rho)\|\leq C|\rho|^{-\ve}\left|\exp(\rho x \overleftarrow{R}_k)\right|, \ |\rho|>T^{-1}, x\in[|\rho|^{-1},T],
$$
where $C$ do not depend on $x,\rho$.

The estimate for $F_k(x,\rho)$ can be obtained in a similar way using the second part of Lemma 3.2 and successive approximation method applied directly to equation (3.4).
$\hfil\Box$

\medskip
{\bf Corollary 3.1.} For any fixed $x\in (0,\infty)$ and $\rho\to\infty$, $\rho\in\overline S_\nu$ the following asymptotics hold:
$$
T_k(x,\rho)=T^0_k(x,\rho)+O\left(\rho^{-\ve}\exp\left(\rho x \overleftarrow{R}_k\right)\right),\ \ve\in(0,1), $$$$
F_k(x,\rho)=F^0_k(x,\rho)+O\left(\rho^{-1}\exp\left(\rho x \overrightarrow{R}_k\right)\right).
$$

\medskip
We complete this section with the following result concerning a behavior of  the fundamental tensors for $\rho\to 0$.
Consider the tensors $\tilde F_k(x,\rho):=\rho^{-\overrightarrow{\mu}_k}F_k(x,\rho)$, $\tilde T_k(x,\rho):=\rho^{-\overleftarrow{\mu}_k}T_k(x,\rho)$.

\medskip
{\bf Lemma 3.4.} For all $k=\overline{1,n}$ $\tilde F_k(x,\cdot)$, $\tilde T_k(x,\cdot)$ admit the continuous extensions to $\overline S_\nu$.

\medskip
{\bf Proof.}
First we notice that $\tilde F^0_k(x,\rho):=\rho^{-\overrightarrow{\mu}_k}F^0_k(x,\rho)$, $\tilde T^0_k(x,\rho):=\rho^{-\overleftarrow{\mu}_k}T^0_k(x,\rho)$ are continuous w.r.t. $\rho$ in $\overline S_\nu$.
Indeed, $$\tilde T^0_k(x,\rho)=\left(x^{\mu_k}\hat c_k(\rho x)\right)\wedge\dots\wedge\left(x^{\mu_n}\hat c_n(\rho x)\right)$$ is an entire function of $\rho$. Further, we write $\tilde F_k(x,\rho)=\tilde \Psi^0_1(x,\rho)\wedge\dots\wedge\tilde \Psi^0_k(x,\rho)$ (where $\tilde \Psi^0_k(x,\rho):=\rho^{-\mu_k}\Psi^0_k(x,\rho)$) and
$$\tilde\Psi^0_j(x,\rho)=\sum\limits_{s=j}^n l_{sj}\rho^{\mu_s-\mu_j}x^{\mu_s}\hat c_s(\rho x)$$
with constant lower-triangular matrix $l$ (see Section 2). Since for $s>j$ one has $\mbox{Re}(\mu_s-\mu_j)>0$, we can conclude that all $\tilde\Psi^0_k(x,\rho)$ and $\tilde F^0_k(x,\rho)$, $k=\overline{1,N}$ are also continuous w.r.t. $\rho\in\overline S_\nu$.

Now we observe that $\tilde F_k(x,\cdot)$, $\tilde T_k(x,\cdot)$ solve the following slightly modified versions of (3.3), (3.4):
$$
\tilde T_k(x,\rho)=\tilde T_k^0(x,\rho)+\int\limits_0^x G_{n-k+1}(x,t,\rho)\left(q^{(n-k+1)}(t)\tilde T_k(t,\rho)\right)dt,
$$
$$
\tilde F_k(x,\rho)=\tilde F_k^0(x,\rho)-\int\limits_x^\infty G_{k}(x,t,\rho)\left(q^{(k)}(t)\tilde F_k(t,\rho)\right)dt.
$$
Moreover, we can write $$
G_m(x,t,\rho)f=
\sum\limits_{\alpha\in\mathcal{A}_m}(-1)^{\sigma_\alpha}\left|f\wedge\tilde\Psi^0_{\alpha'}(t,\rho)\right|\tilde\Psi^0_{\alpha}(x,\rho),
$$
and notice that the Green's operator $G_m(x,t,\rho)$ is also continuous w.r.t. $\rho\in\overline S_\nu$. Thus we can repeat the arguments used in the proof of Theorem 3.1, for  $\rho\in\overline S_\nu$.
$\hfil\Box$

\bigskip
{\bf 4. Weyl-type solutions.}

As in the previous section, we assume $\rho\in S_\nu$ for some (arbitrary) $\nu\in\overline{1,N}$.

\medskip
{\bf Definition 4.1.} Let $k\in\overline{1,N}$ and $\rho\in S_\nu$ be fixed. Function $y(x)$, $x\in(0,\infty)$ is called {\it $k$-th Weyl-type solution} if it satisfy (1.1) and the following asymptotics hold:
$$
y(x)=O(x^{\mu_k}), x\to 0, \ y(x)=\exp(\rho R_k x)(\mathfrak{f}_k+o(1)), x\to\infty.
$$

\medskip
Below in this section we show that the fundamental tensors $F_k(x,\rho)$ constructed in previous section can actually be represented as the wedge products of Weyl-type solutions. Our first step is to show that $F_k(x,\rho)$ and $T_k(x,\rho)$ are decomposable.

\medskip
{\bf Lemma 4.1.} For any $\rho\in\overline S_\nu\setminus\{0\}$ there exist unique sets of absolutely continuous w.r.t. $x\in(0,\infty)$ functions $\{v_1(x,\rho),\dots, v_n(x,\rho)\}$, $\{w_1(x,\rho),\dots,w_n(x,\rho)\}$ such that:
\begin{itemize}
\item $\{v_1,\dots, v_n\}$ are orthogonal and $\{w_1,\dots,w_n\}$ are orthogonal;\\
\item $T_k=w_k\wedge\dots\wedge w_n$, $F_k=v_1\wedge\dots\wedge v_k$; \\
\item the following asymptotics hold: $$v_k=\exp(\rho R_k x)(\mathfrak{f}_k+o(1)), \ x\to\infty, w_k=(\rho x)^{\mu_k}(\mathfrak{g}_k+o(1)), \ x\to 0,$$
where $\mathfrak{g}_n=\mathfrak{h}_n, \ \mathfrak{g}_k-\mathfrak{h}_k\in span\{\mathfrak{g}_j\}_{j>k}$, and vectors $\{\mathfrak{g}_k\}_{k=1}^n$ are orthogonal; \\
\item the following relations hold: $$(w'_k-Q(x,\rho)w_k)\wedge T_{k+1}=0, F_{k-1}\wedge (v'_k-Q(x,\rho)v_k)=0;$$
\item $v_s\wedge F_k=0$, $w_k\wedge T_s=0$ if $s\leq k$.
\end{itemize}

\medskip
{\bf Proof.} Assertion of the lemma is similar to Lemma 7.2 {\cite{BeDT}} and our proof is based in general on the similar ideas. We consider in more detail construction of $\{w_1,\dots,w_n\}$ ,  $\{v_1,\dots,v_n\}$ can be constructed in the same way.

We use an induction in $k=n,n-1,\dots,1$. The first step is obvious, since we set $w_n:=T_n$. Now suppose we have already constructed the vectors $\{w_{k+1},\dots,w_n\}$ having the properties specified in the lemma and such that for $j>k$ $T_j=w_j\wedge\dots\wedge w_n$. Then we apply Lemma 7.1 {\cite{BeDT}} which asserts that $w_j\wedge T_k=0$ is a sufficient condition for solvability of the system $w\wedge T_{k+1}=T_k$, while the additional requirement $(w, w_j)=0$, $j=\overline{k+1,n}$ makes the solution unique. Let us denote it as $w_k$.

By differentiating  the relation $w_k\wedge T_{k+1}=T_k$ and taking into account that $T_m$ satisfies the equation $T'_m=Q^{(n+1-m)} T_m$ we obtain the relation $(w'_k-Q(x,\rho)w_k)\wedge T_{k+1}=0$.

In order to evaluate an asymptotical behavior of $w_k(x,\rho)$ as $x\to 0$ (while $\rho\in \overline S_\nu\setminus\{0\}$ is arbitrarily fixed) we use the following representation:
$$
w_k(x,\rho)=\sum\limits_{j>k}\beta_j(x,\rho)w_j(x,\rho)+\sum\limits_{j\leq k}\gamma_j(x,\rho)C_j(x,\rho). \eqno(4.1)
$$
Then the relation $w_k\wedge T_{k+1}=T_k$ becomes:
$$
\sum\limits_{j\leq k} \gamma_j C_j\wedge T_{k+1}=T_k,
$$
that yields:
$$
\gamma_j \left|C_1\wedge\dots\wedge C_k\wedge T_{k+1}\right|=$$$$
(-1)^{\chi_j}\left|C_1\wedge\dots\wedge C_{j-1}\wedge C_{j+1}\wedge\dots \wedge C_k\wedge T_k\right|, \ \chi_j\in\{0,1\}, \chi_k=0. \eqno(4.2)
$$
Using the asymptotics of $C_j$ and $T_j$ as $x\to 0$ we obtain from (4.2):
$$
\gamma_k=1+o(1), \ \gamma_{j}=o(\left|(\rho x)^{\mu_k-\mu_j}\right|), \ j<k
$$
and therefore:
$$
w_k^1(x,\rho):=\sum\limits_{j\leq k} \gamma_j C_j=C_k(x,\rho)+o(\left|(\rho x)^{\mu_k}\right|)=(\rho x)^{\mu_k}(\mathfrak{h}_k+o(1)). \eqno(4.3)
$$
Now we consider the relation $(w_k,w_j)=0$, $j=\overline{k+1,n}$ that can be written as follows:
$$
\beta_j(w_j,w_j)+(w^1_k,w_j)=0.
$$
In view of (4.3) and the induction assumption this yields:
$$
\beta_j(x,\rho)=(\beta^0_j+o(1)) (\rho x)^{\mu_k-\mu_j}, \ \beta^0_j(\mathfrak{g}_j,\mathfrak{g}_j)+(\mathfrak{h}_k,\mathfrak{g}_j)=0. \eqno(4.4)
$$
Substituting the obtained representations (4.3), (4.4) into (4.1) we obtain:
$$
w_k(x,\rho)=(\rho x)^{\mu_k}(\mathfrak{g}_k+o(1)), \ \mathfrak{g}_k=\mathfrak{h}_k+\sum\limits_{j>k}\beta^0_j\mathfrak{g}_j. \eqno(4.5)
$$
Finally, for any $j>k$ one has: $(\mathfrak{g}_k,\mathfrak{g}_j)=(\mathfrak{h}_k,\mathfrak{g}_j)+\beta^0_j(\mathfrak{g}_j,\mathfrak{g}_j)=0$ by virtue of (4.4).

In order to complete the induction step we should proof that $w_k\wedge T_j\equiv 0$ for any $j\leq k$.  Denote $w_k\wedge T_j\equiv 0=:y$, we find that $y(x)$ solves the following equation:
$$
y'=Q^{(n-j+2)}_0(x,\rho)y+q^{(n-j+2)}(x)y
$$
and admits the estimate:
$$
\|y(x)\|\leq C \left|\left(\rho x\right)^{\overleftarrow{\mu_j}+\mu_k}\right|, \ |\rho x|\leq 1. \eqno(4.6)
$$
This yields the representation:
$$
y(x)=\int\limits_{x_0}^x G_{n-j+2}(x,t,\rho)\left(q^{(n-j+2)}(t)y(t)\right)dt+\sum\limits_{\alpha\in\mathcal{A}_{n-j+2}}A_\alpha \Psi^0_\alpha(x,\rho),  \eqno(4.7)
$$
that is true for any $x_0>0$ provided the constants $A_\alpha$ are such that:
$$
\sum\limits_{\alpha\in\mathcal{A}_{n-j+2}}A_\alpha\Psi^0_\alpha(x_0,\rho)=y(x_0). \eqno(4.8)
$$
From (4.8) one gets:
$$
|A_\alpha|=\left|\left|y(x_0)\wedge\Psi^0_{\alpha'}(x_0,\rho)\right|\right|.
$$
In view of estimates obtained above for $w_k$ and $T_j$ this yields for $|\rho x_0|<1$:
$$
|A_\alpha|\leq C \left|\left(\rho x_0\right)^{\overleftarrow{\mu_j}+\mu_k-\mu_\alpha}\right|.
$$
Taking the limit as $x_0\to 0$ and taking into account that for any $\alpha\in\mathcal{A}_{n-j+2}$ one has
$\mbox{Re}\left(\overleftarrow{\mu_j}+\mu_k-\mu_\alpha\right)=\mbox{Re}\left(\overleftarrow{\mu}_{j-1}-\mu_\alpha+\mu_k-\mu_{j-1}\right)
\geq\mbox{Re}(\mu_k-\mu_{j-1})>0$
we find that all $A_\alpha\to 0$. Thus we arrive at the following homogeneous Volterra equation with respect to $y$:
$$ y(x)=\int\limits_{0}^x G_{n-j+2}(x,t,\rho)\left(q^{(n-j+2)}(t)y(t)\right)dt. \eqno(4.9) $$
From (4.9) and  estimate (4.6)
we obtain successively for $m=1,2,\dots$:
$$
\|y(x)\|\leq \frac{C}{m!}\cdot \left|\left(\rho x\right)^{\overleftarrow{\mu_j}+\mu_k}\right|\cdot\left(\int\limits_0^x \|q(t)\|dt\right)^m,\ |\rho x|\leq 1.
$$
This means that $y(x)=0$ at least for $x\in(0,|\rho^{-1}|)$. Taking into account that $y(x)$ solves a homogeneous first order system of ODE, we conclude that $y(x)\equiv 0$.
$\hfil\Box$

\medskip
Decomposability  of the fundamental tensors established above is an important fact for our further considerations. But this is only a preliminary step, since
the functions $\{v_k\}, \{w_k\}$ {\it do not satisfy } (1.1). Moreover, one can notice that they are not holomorphic with respect to $\rho$. Our next (and key) step consists in passing from $\{v_k\}$ to the Weyl-type solutions.

Let us define $\Delta_k(\rho)=\left|F_{k-1}(x,\rho)\wedge T_k(x,\rho)\right|$. Clear that $\Delta_k$ do depend on $x$.

\medskip
{\bf Theorem 4.1.} For any $\rho\in\overline S_\nu\setminus\{0\}$ such that $\Delta_k(\rho)\neq 0$ there exists a unique function $\psi_k(x,\rho)$ such that:
\begin{itemize}
\item $F_{k-1}\wedge \psi_k =F_k$, $\psi_k\wedge T_k=0$;\\
\item $\psi'_k=Q(x,\rho)\psi_k$ (i.e. $\psi_k$ solves system (1.1));\\
\item if $\rho\in S_\nu$ then the following asymptotics hold:
$$
\psi_k(x,\rho)=\exp(\rho R_k x)(\mathfrak{f}_k+o(1)), x\to\infty, \ \psi_k(x,\rho)=O((\rho x)^{\mu_k}), x\to 0.
$$
\end{itemize}

\medskip
{\bf Proof.} In view of Lemma 4.1 the relations $F_{k-1}\wedge \psi_k =F_k$, $\psi_k\wedge T_k=0$ are equivalent to:
$$
\psi_k-v_k=\sum\limits_{j< k}\beta_j v_j, \ \psi_k=\sum\limits_{j\geq k}\gamma_j w_j,
$$
that is actually some system of linear algebraic equations w.r.t. $\{\beta_j\}_{j<k}, \{\gamma_j\}_{j\geq k}$ while condition $F_{k-1}\wedge T_k\neq 0$ is actually the condition of its unique solvability. Thus, the condition $\Delta_k(\rho)\neq 0$ guarantees the existence and uniqueness of $\psi_k(x,\rho)$ such that the relations $F_{k-1}\wedge \psi_k =F_k$, $\psi_k\wedge T_k=0$ hold.

To evaluate the asymptotics of $\psi_k(x,\rho)$ as $x\to 0$ we use the representation:
$$
\psi_k(x,\rho)=\sum\limits_{j\geq k}\gamma_j(x,\rho)w_j(x,\rho).
$$
The relation $F_{k-1}\wedge \psi_k=F_k$ becomes:
$$
\sum\limits_{j\geq k}\gamma_j F_{k-1}\wedge w_j=F_k,
$$
that yields:
$$
\gamma_j\Delta_k=(-1)^{\chi_j}\left|F_k \wedge w_k \wedge\dots w_{j-1}\wedge w_{j+1}\wedge\dots \wedge w_n \right|, \ \chi_j\in\{0,1\}.
$$
Using the estimates for $F_k$ from Theorem 3.1 and asymptotics for $w_j$ from Lemma 4.1 we obtain $\gamma_j=O\left((\rho x)^{\mu_k-\mu_j}\right)$ for $|\rho x|\leq 1$ and therefore $\psi_k(x,\rho)=O\left((\rho x)^{\mu_k}\right)$ as $x\to 0$.

Now we consider $\psi_k(x,\rho)$ for $x\to\infty$. Here we use the representation:
$$
\psi_k(x,\rho)=\sum\limits_{j=1}^{k-1}\beta_j(x,\rho)v_j(x,\rho)+v_k(x,\rho). \eqno(4.10)
$$
For the coefficients $\beta_j$ we have the following equation:
$$
\sum\limits_{j=1}^{k-1}\beta_j v_j \wedge T_k=0.
$$
As above, we obtain:
$$
\beta_j \Delta_k=(-1)^{\chi_j} \left| v_1\wedge\dots v_{j-1}\wedge v_{j+1}\wedge\dots\wedge v_{k-1}\wedge v_k \wedge T_k\right|, \ \chi_j\in\{0,1\}.
$$
Using the asymptotics of $v_1, v_2,\dots$ from Lemma 4.1 and the estimate for $T_k$ from Theorem 3.1 one gets:
$$
\beta_j \Delta_k=(-1)^{\chi_j}\delta_j+o(\exp(\rho x(R_k-R_j))), \ \delta_j:=
\left| F^0_{k,j} \wedge T_k\right|, $$
$$ \ F^0_{k,j}:=E_1\wedge\dots E_{j-1}\wedge E_{j+1}\wedge\dots\wedge E_{k-1}\wedge E_k
$$
while for $\delta_j$ one has the representation $\delta_j=\delta_j^0+\hat\delta_j$, where:
$$
\delta_j^0=\left| F^0_{k,j} \wedge T^0_k\right|
$$
does not depend on $x,\rho$ and
$$
\hat\delta_j(x,\rho)=\left| F^0_{k,j} \wedge \hat T_k\right|.
$$
Equation (3.3) yields the following representation for $\hat\delta_j$:
$$
\hat\delta_j=\int\limits_0^x \left|F^0_{k,j}(x,\rho)\wedge\left(G_{n-k+1}(x,t,\rho)\left(q^{(n-k+1)}(t)T_k(t,\rho)\right)\right)
dt\right|.
$$
Now, let us notice that for any $\alpha\in\mathcal{A}_{n-k+1}$ and any system $y_1,\dots, y_n$ of solutions of unperturbed equation (2.1) the Wronskian $\left|y_{\alpha}\wedge F^0_{k,j}\right|$ does not depend on $x$. Thus, from the representation of the Green's operator $G_m(x,t,\rho)$ it follows that:
$$
\hat\delta_j(x,\rho)=\hat\delta_{j0}(\rho)+\hat\delta_{j1}(x,\rho),
$$
where:
$$
\hat\delta_{j0}(\rho)=(-1)^{\chi_\alpha}\int\limits_0^{|\rho^{-1}|} \left|\left(q^{(n-k+1)}(t)T_k(t,\rho)\right)\wedge E_{\alpha'}(t,\rho)\right|\cdot\left|E_{\alpha}(x,\rho)\wedge F^0_{k,j}(x,\rho)\right| dt,\ \chi_\alpha\in\{0,1\},
$$
does not depend on $x$ and:
$$
\hat\delta_{j1}(x,\rho)=(-1)^{\chi_\alpha}\int\limits_{|\rho^{-1}|}^x \left|\left(q^{(n-k+1)}(t)T_k(t,\rho)\right)\wedge E_{\alpha'}(t,\rho)\right|\cdot\left|E_{\alpha}(x,\rho)\wedge F^0_{k,j}(x,\rho)\right| dt,\ \chi_\alpha\in\{0,1\}.
$$
Here $\alpha=(j,k+1,\dots,n)$. Since $\|E_{\alpha}(x,\rho)\wedge F^0_{k,j}(x,\rho)\|=1$, we obtain finally:
$$
\left|\hat\delta_{j1}(x,\rho)\right|=\int\limits_{|\rho|^{-1}}^x \|\left(q^{(n-k+1)}(t)T_k(t,\rho)\right)\wedge E_{\alpha'}(t,\rho)\| dt.
$$
Using the estimates for $T_k$ and $E_{\alpha'}$ we obtain:
$$
|\hat\delta_{j1}(x,\rho)|\leq C \int\limits_{|\rho|^{-1}}^x \|q(t)\|\cdot |\exp(\rho t \overleftarrow{R}_k-\rho t R_\alpha)| dt.
$$
Since $\overleftarrow{R}_k-R_\alpha=R_k-R_j$, we can write:
$$
|\hat\delta_{j1}(x,\rho)\exp(\rho x(R_j-R_k))|\leq C \int\limits_{|\rho|^{-1}}^x \|q(t)\|\cdot |\exp(\rho (x-t)(R_j-R_k)| dt.
$$
For any fixed $\rho\in S_n$ and $j<k$ we have $\mbox{Re}\rho R_j<\mbox{Re}\rho R_k$, therefore the right-hand side of the last estimate tends to $0$ as $x\to\infty$. Thus we obtain $\hat\delta_{j1}=o(\exp(\rho x(R_k-R_j)))$ and $\hat \delta_j= o(\exp(\rho x(R_k-R_j)))$. This yields $\beta_j=\beta^0_j + o(\exp(\rho x(R_k-R_j)))$ and taking again into account that $\mbox{Re}\rho R_j<\mbox{Re}\rho R_k$ we obtain $\beta_j v_j=o(\exp(\rho x R_k))$. Substituting this to (4.10) and using Lemma 4.1 completes the proof. $\hfil\Box$

\medskip
{\bf Lemma 4.2.} Let $\rho\in S_\nu$ be such that $\Delta_k(\rho)\neq 0$. Then any solution $y(x)$ of (1.1) satisfying the conditions:
$$
y(x)=O((\rho x)^{\mu_k}), x\to 0, \ y(x)=\exp(\rho R_k x)(\mathfrak{f}_k+o(1)), x\to\infty
$$
necessarily coincides with $\psi_k(x,\rho)$.

\medskip
{\bf Proof.} We consider the following tensors: $Y_-(x):=(y(x)-\psi_k(x,\rho))\wedge T_k(x,\rho)$ and $Y_+(x):=F_{k-1}(x,\rho)\wedge(y(x)-\psi_k(x,\rho))$. Then $Y_-$ satisfies the equation:
$$
Y'_-=Q_0^{(n-k+2)}(x,\rho)Y_-+q^{(n-k+2)}(x)Y_-
$$
and the estimate:
$$
\|Y_-(x)\|\leq C \left|\left(\rho x\right)^{\overleftarrow{\mu_k}+\mu_k}\right|, \ |\rho x|\leq 1.
$$
Proceeding as in the proof of Lemma 4.1 we conclude that $Y_-(x)$ satisfies the following homogeneous Volterra equation:
$$
Y_-(x)=\int\limits_{0}^x G_{n-k+2}(x,t,\rho)\left(q^{(n-k+2)}(t)Y_-(t)\right)dt
$$
that yields $Y_-=0$.

Similarly for $Y_+$ we have:
$$
Y'_+=Q_0^{(k)}(x,\rho)Y_+ +q^{(k)}(x)Y_+, \ Y_+(x)=o\left(\exp\left(\rho x \overrightarrow{R}_k\right)\right), \ x\to\infty.
$$
This yields:
$$
Y_+(x)=\int\limits_{x_0}^x G_{k}(x,t,\rho)\left(q^{(k)}(t)Y_+(t)\right)dt+\sum\limits_{\alpha\in\mathcal{A}_{k}}A_\alpha \Psi^0_\alpha(x,\rho),
$$
$$
|A_\alpha|=\left|\left|Y_+(x_0)\wedge\Psi^0_{\alpha'}(x_0,\rho)\right|\right|.
$$
The asymptotics of $Y_+$ implies $A_\alpha=o(1)$ as $x_0\to\infty$ for any $\alpha\in\mathcal{A}_k$. Thus, taking the limit as $x_0\to\infty$ we arrive at the following homogeneous Volterra equation:
$$
Y_+(x)=-\int\limits_{x}^\infty G_{k}(x,t,\rho)\left(q^{(k)}(t)Y_+(t)\right)dt.
$$
Repeating again the arguments used in the proof of Theorem 4.1 we obtain $Y_+=0$.

Thus we have $F_{k-1}\wedge (y-\psi_k)=0$ and $(y-\psi_k)\wedge T_k=0$. Since $F_{k-1}\wedge T_k \neq 0$, this means $y-\psi_k=0$.
$\hfil\Box$

\medskip
Lemma 4.2 shows that for $\rho\in S_\nu$ such that $\Delta_k(\rho)\neq 0$ the $k$-th Weyl type solution exists, is unique and coincides with $\psi_k(x,\rho)$.

The rest of this section is devoted to investigation of the properties of Weyl-type solutions as the functions of spectral parameter $\rho$.

\medskip
{\bf Theorem 4.3.} For $\rho\to\infty$, $\rho\in\overline S_\nu$ the functions $\Delta_k(\rho)$ admit the asymptotics:
$$
\Delta_k(\rho)=\Delta_k^0+O\left(\rho^{-\ve}\right), \ \ve\in(0,1).
$$
For any fixed $x$ and $\rho\to\infty$, $\rho\in\overline S_\nu$ the following asymptotics hold:
$$
\psi_k(x,\rho)=\sum\limits_{j=1}^k \gamma^0_{jk} \exp(\rho x R_j) \mathfrak{f}_j + O\left(\rho^{-\ve}\exp(\rho x R_k)\right)
$$
with constants $\gamma^0_{jk}$, $\gamma_{kk}^0=1$ that do not depend on $q(\cdot)$.

\medskip
{\bf Proof.} Using the representation
$
\Delta_k(\rho)=\left|F_{k-1}(x,\rho)\wedge T_k(x,\rho)\right|,
$
Theorem 3.1 and Corollary 3.1 we obtain the desired asymptotics for $\Delta_k(\rho)$ via direct calculation.

Now we use the induction in $k$ to prove the asymptotics for $\psi_k(x,\rho)$. For $k=1$ the asymtotics follows directly  from Corollary 3.1. Now we suppose that for some $k>1$ all the functions $\psi_1(x,\rho),\dots,\psi_{k-1}(x,\rho)$ are already shown to have the asymptotics specified in the theorem. The induction step consists in obtaining the asymptotics for $\psi_k(x,\rho)$.

First we notice that for any fixed $x$ and $\rho\to\infty$, $\rho\in\overline S_\nu$ we have:
$$
\left|\psi_1\wedge\dots\wedge\psi_{k-1}\wedge E_k\wedge\dots\wedge E_n\right|=
\left| E_1\wedge\dots\wedge E_n+O\left(\rho^{-\ve}\right)\right|=
\det\Pi_\nu+O\left(\rho^{-\ve}\right).
$$
Thus for any sufficiently large $\rho$ the vectors $\psi_1,\dots,\psi_{k-1}, E_k,\dots, E_n$ form a basis in $\mathbb{C}^n$ and we can use the expansion:
$$
\psi_{k}(x,\rho)=\sum\limits_{j<k} \beta_{jk}(x,\rho) \psi_j(x,\rho)+\sum\limits_{j\geq k} \beta_{jk}(x,\rho) E_j(x,\rho). \eqno(4.11)
$$
Then the relation $F_{k-1}\wedge\psi_k=F_k$ becomes:
$$
\sum\limits_{j\geq k} \beta_{jk} F_{k-1}\wedge E_j=F_k
$$
and we find:
$$
\beta_{kk}\left|F_{k-1}\wedge E_k \wedge\dots\wedge E_n\right|=\left|F_{k}\wedge E_{k+1} \wedge\dots\wedge E_n\right|,
$$
$$
|\beta_{jk}|\cdot\left\|F_{k-1}\wedge E_k \wedge\dots\wedge E_n\right\|=
\left\|F_k \wedge E_k\wedge\dots\wedge E_{j-1}\wedge E_{j+1}\wedge\dots E_n\right\|, \ j>k .
$$
Using Corollary 3.1 and taking into account that $\left|F^0_{k-1}\wedge E_k \wedge\dots\wedge E_n\right|=\left|F^0_{k}\wedge E_{k+1} \wedge\dots\wedge E_n\right|$, $F^0_k \wedge E_k=0$ we calculate:
$$
\beta_{kk}(x,\rho)=1+O\left(\rho^{-1}\right), \ \beta_{jk}(x,\rho)=O\left(\rho^{-1}\exp(\rho x (R_k-R_j))\right), $$$$
\psi^1_k(x,\rho):=\sum\limits_{j\geq k} \beta_{jk}(x,\rho) E_j(x,\rho)=E_k(x,\rho)+O\left(\rho^{-1}\exp(\rho x R_k)\right). \eqno(4.12)
$$
We rewrite (4.11) as follows:
$$
\psi_{k}(x,\rho)=\sum\limits_{j<k} \beta_{jk}(x,\rho) \psi_j(x,\rho)+\psi^1_k(x,\rho). \eqno(4.13)
$$
Then the relation $\psi_k\wedge T_k=0$ becomes:
$$
\sum\limits_{j<k}\beta_{jk} \psi_j \wedge T_k= -\psi^1_k \wedge T_k
$$
and we calculate:
$$
(-1)^{\chi_{jk}}\beta_{jk} \Delta_k=\left| \tilde F_{k-1,j} \wedge \psi^1_k \wedge T_k\right|, \ \chi_{jk}\in\{0,1\}. \eqno(4.14)
$$
Here we denote (for $j\leq m$) $\tilde F_{m,j}:=\psi_1 \wedge\dots\psi_{j-1}\wedge\psi_{j+1}\wedge\dots\wedge\psi_m$ and take into account that $\left|F_{k-1}\wedge T_k\right|=\Delta_k$. Now we notice that the induction assumption implies, in particular, $\psi_j=\Psi^0_j+O\left(\rho^{-\ve}\exp(\rho x R_j)\right)$, $j<k$ (since coefficients in the asymptotics do not depend on $q(\cdot)$ and should be the same for the unperturbed equation). Therefore we have:
$$
\tilde F_{k-1,j}(x,\rho)=\tilde F^0_{k-1,j}(x,\rho)+O\left(\rho^{-\ve}\exp(\rho x(\overrightarrow{R}_{k-1}-R_j))\right),
$$
$$
\tilde F_{k-1,j}\wedge\psi^1_k=\tilde F^0_{k-1,j}\wedge E_k+O\left(\rho^{-\ve}\exp(\rho x(\overrightarrow{R}_{k}-R_j))\right),
$$
where $\tilde F^0_{m,j}:=\Psi^0_1 \wedge\dots\Psi^0_{j-1}\wedge\Psi^0_{j+1}\wedge\dots\wedge\Psi^0_m$.
Thus, we can use Corollary 3.1 and Theorem 3.1 to obtain from (4.14) the following asymptotics:
$$
(-1)^{\chi_{jk}}\beta_{jk} \Delta_k=\left|\tilde F^0_{k-1,j}\wedge E_k\wedge T^0_k \right|+
O\left(\rho^{-\ve}\exp(\rho x(R_k-R_j))\right). \eqno(4.15)
$$
We notice again that the Wroskian $\left|\tilde F^0_{k-1,j}\wedge E_k\wedge T^0_k \right|$ of the solutions of (2.1) does not depend on $x,\rho$. Then from (4.15) it follows that:
$$
\beta_{jk}(x,\rho)=\beta^0_{jk}+O\left(\rho^{-\ve}\exp(\rho x(R_k-R_j))\right),
$$
where $\beta^0_{jk}$ obviously do not depend on $q(\cdot)$. Now we return to (4.12), (4.13) and use again the induction assumption. This completes the proof.
$\hfil\Box$

\medskip
Since the fundamental tensors $T_k(x,\rho)$, $F_k(x,\rho)$ $k=\overline{1,n}$ are analytical function of $\rho\in S_\nu$ and continuous in $\overline S_\nu\setminus\{0\}$, we can use our construction of the Weyl-type solutions to conclude that $\psi_k(x,\cdot)$ is holomorphic in $S_\nu\setminus Z_{\nu k}$ and continuous in $\overline S_\nu\setminus \left(Z_{\nu k}\cup\{0\}\right)$, where $Z_{\nu k}$ is a set of zeros of $\Delta_k(\cdot)$ in $\overline S_\nu\setminus\{0\}$.

The observations of Lemma 3.4 above allow us to evaluate a behavior of $\psi_k(x,\rho)$ for $\rho\to 0$. First we notice that $\overrightarrow{\mu}_{k-1}+\overleftarrow{\mu}_k=\sum\limits_{j=1}^n \mu_j=0$ and therefore $\Delta_k(\rho)=F_{k-1}(x,\rho)\wedge T_k(x,\rho)=\tilde F_{k-1}(x,\rho)\wedge \tilde T_k(x,\rho)$ are continuous w.r.t. $\rho\in\overline S_\nu$.

\medskip
{\bf Theorem 4.4.} If for all $k=\overline{1,n}$ $\Delta_k(0)\neq 0$ then all the functions $\rho^{-\mu_k}\psi_k(x,\rho)$ admit the continuous w.r.t. $\rho$ extensions to $\overline S_\nu\cap\{|\rho|\leq\delta\}$ for some positive $\delta$.

\medskip
{\bf Proof.} Our consideration uses induction in $k$. Suppose that $\tilde \psi_j(x,\rho):=\rho^{-\mu_j}\psi_j(x,\rho)$, $j=\overline{1,k-1}$ are already shown to be continuous in $\rho\in\overline S_\nu\cap\{|\rho|\leq \delta\}$ for some $\delta>0$. Then the induction step consists in proving that $\tilde\psi_k(x,\rho)$ is also continuous in $\rho\in\overline S_\nu\cap\{|\rho|\leq \delta\}$.

Further we assume that $x>0$ is arbitrary fixed and omit it in all the argument's lists.

At this point we have:
$$
\tilde F_{k}(\rho)=\tilde\psi_1(\rho)\wedge\dots\wedge\tilde\psi_{k}(\rho), \rho\in\overline S_\nu\cap\{|\rho|\leq \delta\}\setminus\{0\}, \eqno(4.16)
$$
while the induction assumption yields:
$$
\tilde F_{k-1}(\rho)=\tilde\psi_1(\rho)\wedge\dots\wedge\tilde\psi_{k-1}(\rho), \rho\in\overline S_\nu\cap\{|\rho|\leq \delta\}.
$$
From (4.16) it follows that $$\tilde\psi_j(\rho)\wedge\tilde F_k(\rho)=0, \ j=\overline{1,k-1}, \rho\in\overline S_\nu\cap\{|\rho|\leq \delta\}\setminus\{0\} \eqno(4.17).$$
But the left-hand side of (4.17) is continuous in $\rho\in\overline S_\nu\cap\{|\rho|\leq \delta\}$ therefore the relation is true for $\rho=0$ as well. By virtue of Lemma 7.1 {\cite{BeDT}} there exists a unique vector-function $f(\rho)$ that is orthogonal to all $\tilde\psi_j(\rho), j=\overline{1,k-1}$ and such that $\tilde F_{k-1}\wedge f=\tilde F_{k}$ for $\rho\in\overline S_\nu\cap\{|\rho|\leq \delta\}$.

Since $\Delta_{k-1}(0)\neq 0$, it is possible to choose some fixed vectors $\mathfrak{u}_k,\dots,\mathfrak{u}_n$ such that $\tilde\psi_1(\rho)\wedge\dots\wedge\tilde\psi_{k-1}(\rho)\wedge\mathfrak{u}_k\wedge\dots\wedge\mathfrak{u}_n\neq 0$ for all sufficiently small $\rho\in\overline S_\nu$. Therefore for all sufficiently small $\rho\in\overline S_\nu$ $f(\rho)$ admits the following representation:
$$
f(\rho)=\sum\limits_{j<k} \gamma_j(\rho)\tilde\psi_j(\rho)+\sum\limits_{j\geq k}\gamma_j(\rho)\mathfrak{u}_j.
$$
Relations $(f,\tilde\psi_j)=0$, $j=\overline{1,k-1}$, $\tilde F_{k-1}\wedge f=\tilde F_k$ are equivalent to some system of linear equations with respect to $\gamma_j$, $j=\overline{1,n}$ with the coefficients, which are continuous w.r.t. $\rho\in\overline S_\nu$. Condition $\tilde\psi_j\wedge\tilde F_k=0$, $j=\overline{1,k-1}$ allow us to convert this system into some $(k-1)\times(k-1)$ linear system, which is non-degenerate for all sufficiently small  $\rho\in\overline S_\nu$.
 Thus we can conclude that $f(\rho)$ is continuous in $\rho\in\overline S_\nu \cap \{|\rho|\leq\delta\}$ with some $\delta>0$. Since for $\rho \in \overline S_\nu\setminus\{0\}$ we have $\tilde F_{k-1}\wedge (\tilde\psi_k-f)=0$, we find that $\tilde\psi_k$ admits the representation:
$$
\tilde\psi_k(\rho)=f(\rho)-\sum\limits_{j<k} \beta_j(\rho)\tilde\psi_j(\rho), \
\rho\in\overline S_\nu\cap\{|\rho|\leq \delta\}\setminus\{0\}. \eqno(4.18)
$$
Then the relation $\tilde\psi_k \wedge \tilde T_k=0$ becomes:
$$
\sum\limits_{j<k} \beta_j(\rho)\tilde\psi_j(\rho)\wedge\tilde T_k(\rho)=f(\rho)\wedge\tilde T_k(\rho),
$$
that yields
$$
\beta_j(\rho)\Delta_k(\rho)=$$$$(-1)^{\chi_j}\left| f(\rho)\wedge\tilde\psi_1(\rho)\wedge\dots\wedge\tilde\psi_{j-1}(\rho)\wedge\tilde\psi_{j+1}(\rho)\dots\wedge\tilde\psi_{k-1}(\rho)
\wedge\tilde T_k(\rho)\right|, \ \chi_j\in\{0,1\}.
$$
By virtue of the condition $\Delta_k(0)\neq 0$ and of the induction assumption this relation defines the coefficients $\beta_j(\rho)$ as some continuous function of $\rho\in \overline S_\nu\cap \{|\rho|\leq\delta\}$. Thus the right-hand side in (4.18) is a continuous function of  $\rho\in \overline S_\nu\cap \{|\rho|\leq\delta\}$ and the induction step is verified.
$\hfil\Box$

\bigskip
{\bf 5. Inverse scattering problem. Uniqueness result.}

Here we consider the inverse scattering problem for system (1.1). For the sake of simplicity, we restrict our consideration with "solitonless" case. Specifically, we assume that the following condition is satisfied.

\medskip
{\bf Condition $G_0$.} For any sector $S_\nu$, $\nu=\overline{1,N}$ all the functions $\Delta_k(\rho)$, $k=\overline{1,n}$ do not vanish for $\rho\in\overline S_\nu$.

\medskip
For each sector $S_\nu$, $\nu=\overline{1,N}$ we consider the fundamental matrix $$\Psi(x,\rho):=(\psi_1(x,\rho),\dots,\psi_n(x,\rho)),$$ where $\psi_k(x,\rho)$ are the Weyl-type solutions of (1.1) constructed in the previous section. Let $\Sigma_\nu$ be a ray that separates the sectors $S_\nu$ and $S_{\nu+1}$ (we assume that the sectors $S_1, S_2,\dots$ are enumerated in counterclockwise order and $S_{N+1}:=S_1$). Under Condition $G_0$ for any $\rho\in \Sigma_\nu$ there exist the boundary values $\Psi_-(x,\rho):=\lim\limits_{\xi\to\rho, \xi\in S_\nu}\Psi(x,\xi)$ and $\Psi_+(x,\rho):=\lim\limits_{\xi\to\rho, \xi\in S_{\nu+1}}\Psi(x,\xi)$. We define for $\rho\in \Sigma_\nu$ the matrix $v(\rho):=\Psi^{-1}_-(x,\rho)\Psi_+(x,\rho)$. Following the general scheme of {\cite{BeDT}} we treat $v(\cdot)$ as {\it scattering data} and consider the following inverse scattering problem.

\medskip
{\bf Problem IP0.} Given $v(\rho), \rho\in\Sigma\setminus\{0\}$, $A,B$ recover $q(x), x\in(0,\infty)$.

\medskip
In the sequel together with system (1.1) with the matrix coefficient $q(\cdot)$ we consider the system of the same form but with different (in general) coefficient $\tilde q(\cdot)$. We agree that if a symbol $\xi$ denotes some object related to system (1.1) with the coefficient $q(\cdot)$ then $\tilde\xi$ denotes an analogous object related to the system with the coefficient $\tilde q(\cdot)$.

Main result of this section is the following theorem that asserts that under  conditions $R_0$ and $G_0$ specification of the matrix $v(\rho)$, $\rho\in\Sigma\setminus\{0\}$ determines uniquely the coefficient $q(x)$, $x\in(0,\infty)$.

\medskip
{\bf Theorem 5.1.} Let the matrix $A,B$ be such that Condition $R_0$ is satisfied. Let $q(\cdot)$ and $\tilde q(\cdot)$ be such that Condition $G_0$ is satisfied. Then $\tilde v(\rho)=v(\rho)$ for all $\rho\in\Sigma\setminus\{0\}$ implies $\tilde q(x)=q(x)$ for a.e. $x\in(0,\infty)$.

\medskip
{\bf Proof.} Let us consider the following {\it spectral mappings matrix}{\cite{Ybook}}:
$$
P(x,\rho):=\Psi(x,\rho)\tilde\Psi^{-1}(x,\rho).
$$
Under the conditions of theorem for each fixed $x\in(0,\infty)$ $P(x,\cdot)$ is a sectionary-holomorphic function with jump matrices on $\Sigma_\nu$, $\nu=\overline{1,N}$ equal to:
$$
P^{-1}_-(x,\rho)P_+(x,\rho)=\tilde\Psi_-(x,\rho)v(\rho)\tilde v^{-1}(\rho)\tilde\Psi^{-1}_-(x,\rho).
$$
If $\tilde v=v$ then we have $P_+(x,\rho)=P_-(x,\rho)$ for all $\rho\in\Sigma\setminus\{0\}$. Therefore $P(x,\cdot)$ is holomorphic in all $\rho$-plane with possible exception of $\rho=0$. Moreover, by virtue of Theorem 4.4 $P(x,\rho)=\left(\Psi(x,\rho)\rho^{-\mu}\right)\left(\tilde\Psi(x,\rho)\rho^{-\mu}\right)^{-1}$ (where $\mu:=diag(\mu_1,\dots,\mu_n)$) is bounded as $\rho\to 0$. Furthermore, by virtue of Theorem 4.3 we have
$$
\Psi(x,\rho)=\left(\mathfrak{f}\Gamma^0(x,\rho)+O(\rho^{-\ve})\right)\exp(\rho x R),
$$
where $\mathfrak{f}:=(\mathfrak{f}_1,\dots,\mathfrak{f}_n)$, $R:=diag(R_1,\dots,R_n)$, $\Gamma^0_{jk}(x,\rho)=\gamma^0_{jk}\exp(\rho x(R_j-R_k))$ for $j\leq k$ and $\Gamma^0_{jk}(x,\rho)=0$ for $j>k$. Since for $\rho\in\overline S_\nu$ and $j\leq k$ we have $\mbox{Re}(\rho(R_j-R_k))\leq 0$ and $\gamma^0_{kk}=1$, we can estimate $\Gamma^0(x,\rho)=O(1)$ and $\left(\Gamma^0(x,\rho)\right)^{-1}=O(1)$. Therefore
$$
\Psi(x,\rho)=\mathfrak{f}\Gamma^0(x,\rho)\left(I+O(\rho^{-\ve})\right)\exp(\rho x R).
$$
Since the matrix $\Gamma^0(x,\rho)$ does not depend on $q(\cdot)$, this yields $P(x,\rho)=I+O\left(\rho^{-\ve}\right)$ as $\rho\to\infty$. Thus, we can conclude that $P(x,\rho)\equiv I$, $\tilde\Psi(x,\rho)=\Psi(x,\rho)$ and finally $\tilde q(x)=q(x)$ for a.e. $x$.
$\hfil\Box$

\medskip
{\bf Acknowledgement.}
This work was supported by Grant 1.1436.2014K of the Russian Ministry of
Education and Science and by Grants 15-01-04864 and 16-01-00015 of Russian
Foundation for Basic Research.

\noindent Ignatyev, Mikhail\\
Department of Mathematics, Saratov State University, \\
Astrakhanskaya 83, Saratov 410012, Russia, \\
e-mail: mikkieram@gmail.com, ignatievmu@info.sgu.ru

\end{document}